\documentclass[12pt]{amsart}
\usepackage{amssymb}
\usepackage{amsmath}
\usepackage{amsthm}
\usepackage{times}
\usepackage{graphicx}
\usepackage{changepage}
\usepackage{caption}
\usepackage{geometry}
\usepackage[curve,all]{xy}
 \xyoption{curve}
 \xyoption{line}
\xyoption{arc}
\xyoption{color}
\usepackage{array}
\usepackage{enumitem}

\usepackage{color}
\usepackage{amsthm}
\newtheorem{theorem}{Theorem}[section]
\newtheorem{lemma}[theorem]{Lemma}
\newtheorem{corollary}[theorem]{Corollary}

\usepackage{subfig}
 
\theoremstyle{definition}

\newtheorem{example}[theorem]{Example}

\theoremstyle{remark}

\newtheorem{question}[theorem]{Question}

\numberwithin{equation}{section}

\def\deg{{\text{deg}}}

\begin{document}
\title[Quasi-elliptic surfaces]{On the multicanonical systems of 
quasi-elliptic surfaces in characteristic 3}

\author{Toshiyuki Katsura}\thanks{Research of the author is partially supported by
Grant-in-Aid for Scientific Research (B) No. 15H03614.}

\address{Faculty of Science and Engineering, Hosei University,
Koganei-shi, Tokyo 184-8584, Japan}
\email{toshiyuki.katsura.tk@hosei.ac.jp}

\begin{abstract}
We consider the multicanonical systems $\vert mK_S \vert$  of quasi-elliptic surfaces 
with Kodaira dimension $1$ in characteristic 3. 
We show that for any $m \geq 5$ $\vert mK_S \vert$  
gives the structure 
of quasi-elliptic fiber space, and the number $5$ is best possible to give the structure for
any such surfaces.
\end{abstract}

\maketitle

Dedicated to Professor Piotr Pragacz on the occasion of 
his sixtieth birthday

\maketitle

\section{Introduction}
Let k be an algebraically closed field in characteristic $p > 0$, and
let $\varphi : S \longrightarrow B$ be an elliptic surface
with Kodaira dimension $\kappa (S) = 1$ over $k$. 
Let $K_{S}$ be a canonical
divisor of $S$. We consider the multicanonical system $\vert mK_S \vert$.
In Katsura-Ueno \cite{KU} and Katsura \cite{K}, we considered 
the following question:

\begin{question}
Is there a positive integer $M$ such that if $m\geq M$,
then the multicanonical system $\vert mK_S \vert$  gives a structure
of elliptic surface for any elliptic surface $S$ over $k$ with $\kappa (S) = 1$?
\end{question}

For the complex analytic case, Iitaka showed that $M = 86$
and 86 is best possible (cf. Iitaka \cite{I}). Namely, if
$m$ is smaller than 86, then there exists an elliptic surface
$\varphi : S \longrightarrow B$
with $\kappa (S) = 1$ such that  $\vert mK_S \vert$ 
does not give the structure of elliptic surface.
In the case of algebraic
elliptic surfaces, if the characteristic $p = 0$ or $p \geq 3$, 
we showed that $M = 14$ and 14 is best possible 
(Katsura-Ueno \cite{KU} and Katsura \cite{K}). If $p =2$, then we showed
that $M = 12$ and 12 is best possible (cf. Katsura \cite{K}).

In this paper, we treat quasi-elliptic surfaces $\varphi : S \longrightarrow B$
and consider the following similar question.

\begin{question}
Is there a positive integer $M$ such that if $m\geq M$,
then the multicanonical system $\vert mK_S \vert$  gives the structure
of quasi-elliptic surface for any quasi-elliptic surface $S$ over $k$ with $\kappa (S) = 1$?
\end{question}

Note that
quasi-elliptic surfaces exist only in characterisitics $p = 2$ and $3$
(cf. Bombieri-Mumford \cite{BM}). We  show that
if $p = 3$, then we have $M = 5$ and the number $5$ is best possible.
In characteristic $ p = 2$, we have still some difficulties
to determine the best possible number. 

In Section 2, we summarize basic facts on quasi-elliptic surfaces.
In Section 3, we examine the multicanonical system of
quasi-elliptic surfaces in characteristic 3 and show our main theorem.

\section{Some lemmas for quasi-elliptic surfaces}
Let k be an algebraically closed field in characteristic $p > 0$ and let 
$\varphi : S \longrightarrow B$ be a quasi-elliptic surface defined over $k$.
Throughout this paper, we assume that any exceptional curve of the first kind
is not contained in fibers.
Such a surface exists if and only if $p$ is equal to 2 or 3, and
the multiplicities of multiple fibers are all equal to $p$ (cf. Bombieri-Mumford \cite{BM}).
We denote by $g$ the genus of the curve $B$.
In this section, 
we recall some facts on quasi-elliptic surfaces.
We denote by $pF_{i}$ $(i = 1, \ldots,\lambda)$ the multiple fibers. 
Let ${\mathcal T}$ be the torsion part of $R^{1}\varphi_{*}{\mathcal O}_{S}$.
Then, there exists a Cartier divisor ${\bf f}$ on $B$ such that 
$R^{1}\varphi_{*}{\mathcal O}_{S}/{\mathcal T} \cong {\mathcal O}_{B}({\bf f})$.
The canonical divisor formula of $S$ is given by
$$     
K_{S} \sim \varphi^{*}(K_{B} - {\bf f}) + \sum_{i= 1}^{\lambda}a_{i}F_{i},
$$
where $- \deg {\bf f} = \chi({{\mathcal O}_{S}})  + t$ with 
$t =$ the length of 
the torsion part of $R^{1}\varphi_{*}{\mathcal O}_{S}$ and $0 \leq a_{i} \leq p - 1$. 
Here, $\sim$ means linear equivalence.
If $pF_i$ is a tame multiple fiber, then we have $a_i = p -1$.
For details, see Bombieri-Mumford~\cite{BM}.

\begin{lemma}\label{lm:alb}
The Albanese variety ${\rm Alb}(S)$ of $S$ is isomorphic to the Jacobian variety ${\rm J}(B)$
of $B$.
\end{lemma}
\proof{Let $\psi : S \longrightarrow {\rm Alb}(S)$ be the Albanese mapping.
If ${\rm Alb}(S)$ is a point, then by the universality of Albanese variety
we see that the Jacobian variety ${\rm J}(B)$ of $B$ is also a point.
Now, assume ${\rm Alb}(S)$ is not a point.
Since the general fiber of $\varphi$ is a rational curve with one cusp,
the fibers are contracted by $\psi$. Therefore, $\psi (S)$ is a curve. 
We have, by the universality of Albanese variety, a commutative diagram:
$$
\begin{array}{rccc}
\varphi : &  \quad S  & \longrightarrow  & B\\
   & \psi \downarrow &       & \downarrow \\
   & {\rm Alb}(S) & \longrightarrow & {\rm J}(B) \\
   & \quad  \cup  &       & \cup \\
   & \quad  \psi (S)  &   & B.
\end{array}
$$
By this diagram, we have a morphism $\psi (S) \longrightarrow B$.
Therefore, by the Stein factorization theorem, we see that
$\psi (S)$ is isomorphic to $B$. Therefore, by the universality of Jacobian variety,
we conclude ${\rm Alb}(S) \cong {\rm J}(B)$ (see also Katsura-Ueno \cite{KU}, Lemma 3.4).}

We find the following lemma and corollary 
in Lang \cite{L} and Raynaud \cite{R}.
We give here an easy proof for the lemma.
\begin{lemma}\label{lm:chi}
Let $\varphi : S \longrightarrow B$ be a quasi-elliptic surface over a non-singular complete curve
$B$ with genus $g$. Then, we have the inequality
$\chi ({\mathcal O}_{S}) \geq (1 -g)/3$.
\end{lemma}
\proof{By Noether's formula and the self-intersection number $c_1(S)^2 = 0$ of the first Chern class of $S$, we have
$$
12 \chi ({\mathcal O}_{S}) = c_{1}(S)^{2} + c_{2}(S) = 2 - 4q(S) + b_{2}(S).
$$
By Lemma~\ref{lm:alb}, we have $q(S) = g$. Denoting by $\rho (S)$
the Picard number of $S$, we have also $b_{2}(S) \geq \rho (S) \geq 2$. 
Hence, we have $\chi ({\mathcal O}_{S}) \geq (1 -g)/3$.}
\begin{corollary}\label{cor:chi}
{\rm (i)} If $g = 1$, then $\chi ({\mathcal O}_{S}) \geq 0$.

{\rm (ii)} If $g = 0$, then $\chi ({\mathcal O}_{S}) \geq 1$.
\end{corollary}

\section{Multicanonical systems}
In this section, let $k$ be an algebraically closed field of 
characteristic 3.
Let $\varphi: S \longrightarrow B$ be a quasi-elliptic surface defined over $k$.

\begin{example}\label{ex:3} 
In characteristic 3, we consider the quasi-elliptic surface $\varphi : S \longrightarrow {\bf P}^1$ 
which is given by 
a non-singular complete model of the surface defined by
$$
t^{2}(t-1)^{2}z + t^{2}(t -1)^{2}+x^{3} +tz^{3} = 0
$$
Here, $t$ is a parameter of the base curve ${\bf P}^{1}$. 
By Lang \cite{L} p.485, this surface has 
two tame multiple fibers at $t = 0, 1$,
and we have $\chi({\mathcal O}_{S}) = 1$. We denote the two tame multiple
fibers by $3F_0$ and $3F_1$. The canonical divisor $K_S$ is given by
$$
\begin{array}{rl}
K_S & \sim \varphi^* (K_{{\bf P}^{1}} - {\bf f}) + 2F_0 + 2F_1\\
      & \sim - F + 2F_0 + 2F_1.
\end{array}
$$
Here, ${\bf f}$ is a Cartier divisor on ${\bf P}^1$ with $- \deg {\bf f} = \chi({\mathcal O}_S) = 1$
such that 
${\mathcal O}_{{\bf P}^1}({\bf f}) \cong R^{1}\varphi_{*}{\mathcal O}_{S}$,  and $F$ is a general fiber of 
$\varphi : S \longrightarrow {\bf P}^1$. Since we have $F \sim 3F_0 \sim 3F_1$,
we see $4K_S \sim 2F_0 + 2F_1$. Therefore, we have $\dim {\rm H}^0(S, {\mathcal O}_{S}(4K_S)) = 1$, and $\vert 4 K_{S} \vert$ does not give 
the structure of quasi-elliptic surface. If $m  \geq 5$, then we have 
$\dim {\rm H}^0(S, {\mathcal O}_{S}(mK_S)) \geq 2$, and $\vert m K_{S} \vert$ gives 
the structure of quasi-elliptic surface.
\end{example}
We have the following theorem.
\begin{theorem}
Assume that the characteristic $p = 3$. Then, for any quasi-elliptic surface 
$f : S \to B$ with $\kappa (S) = 1$ over $k$ and
for any $ m \geq 5$, the multicanonical system $\vert mK_S \vert$  
gives the unique structure of quasi-elliptic surface, and the number 5 is best possible.
\end{theorem}
\proof{The method of the proof is similar to the one in Iitaka~\cite{I}, 
Katsura-Ueno \cite{KU} and Katsura~\cite{K}. Since the Kodaira dimension
is equal to 1, the structure of quasi-elliptic surface is unique. 
The Kodaira dimension of $S$ is equal to 1 if and only if
$$
(*) \quad 2g -2 + \chi({\mathcal O}_{S}) + t + \sum_{i=1}^{\lambda}(a_{i}/3) > 0.
$$
Therefore, we need to find the least integer $m$ such that
$$
(**)\quad m(2g -2 + \chi({\mathcal O}_{S}) + t) + \sum_{i = 1}^{\lambda}[ma_{i}/3]
\geq 2g + 1
$$
holds under the condition $(*)$. Here, $[r]$ means the integral part of 
a real number $r$. 
We have the following 6 cases:

Case (I) $g \geq 2$

Case (II-1) $g = 1, \chi({\mathcal O}_{S}) + t \geq 1$

Case (II-2) $g = 1, \chi({\mathcal O}_{S}) = 0, t = 0$

Case (III-1) $g = 0, \chi({\mathcal O}_{S}) + t \geq 3$

Case (III-2) $g = 0, \chi({\mathcal O}_{S}) + t = 2$

Case (III-3) $g = 0, \chi({\mathcal O}_{S})= 1, t = 0$

\noindent
We check $(**)$ under the condition $(*)$ for each case.
In Case (I), by Lemma~\ref{lm:chi}, we have $2g - 2 +  \chi({\mathcal O}_{S})\geq 5(g - 1)/3$.
Hence, if $m \geq 3$, $(**)$ holds.
In Case (II-1), if $m \geq 3$,  $(**)$ holds.
In Case (II-2), 
all multiple fibers are tame in this case, and we have at least one multiple fiber by $(*)$.
Since $a_i = 2$, $(**)$ holds for $m \geq 5$.
In Case (III-1), $(**)$ holds for $m \geq 1$.
In Case (III-2), 
since $\chi ({\mathcal O}_{S}) \geq 1$ by Corollary \ref{cor:chi}, we have $t \leq 1$. 
Therefore, the number of wild fibers is less than or equal to $1$.
If there exists at least one tame multiple fiber, then $(**)$ holds for $m \geq 2$. 
If there exist no tame fibers and only one wild fiber, then by $(*)$ 
we have  $a_1 \geq 1$. Therefore,
$(**)$ holds for 
$m \geq 3$.
In Case (III-3), all multiple fibers are tame, and we have $\lambda \geq 2$ by $(*)$. 
Therefore, $(**)$ holds for $m \geq 5$.
The result on the best possible number in characteristic 3 follows 
from Example \ref{ex:3}.
}

In characteristic 2, we can also consider a similar question to the one in characteristic 3.
We have still difficulties to decide the best possible number. For example we need to solve
the following question.

\begin{question}
Does there exist a quasi-elliptic surface over an elliptic curve 
with only one tame multiple fiber  and with $\chi({\mathcal O}_{S}) = 0$ in characteristic 2?
\end{question}

If there don't exist such quasi-elliptic surfaces, then we can show that
in characteristic 2, $(**)$ holds for $m \geq 4$
and that the best possible number is equal to $4$. Namely, we have $M = 4$ in characteristic 2.

\end{document}